\documentclass[aos,preprint]{imsart}

\RequirePackage[OT1]{fontenc}
\RequirePackage{amsthm,amsmath,amsfonts}
\RequirePackage[numbers]{natbib}
\RequirePackage[colorlinks,citecolor=blue,urlcolor=blue]{hyperref}
\RequirePackage{hypernat}
\RequirePackage[pdftex]{graphics}
\RequirePackage{graphicx}
\RequirePackage{float}

\startlocaldefs
\numberwithin{equation}{section}
\theoremstyle{plain}
\newtheorem{theorem}{Theorem}[section]
\newtheorem{lemma}{Lemma}[section]
\newtheorem{proposition}{Proposition}[section]

\newtheorem{example}{Example}[section]
\newtheorem{remark}{Remark}[section]
\newtheorem{assumption}{Assumption}[section]

\startlocaldefs
\endlocaldefs

\begin{document}

\begin{frontmatter}

\title{Surprising Asymptotic Conical Structure in Critical Sample Eigen-Directions}
\runtitle{Conical Structure of Sample Eigen-Directions}


\begin{aug}
\author{\fnms{Dan} \snm{Shen}\thanksref{m1,t1}\ead[label=e1]{dshen@email.unc.edu}},
\author{\fnms{Haipeng} \snm{Shen}\thanksref{t2}\ead[label=e2]{haipeng@email.unc.edu}},
\author{\fnms{Hongtu} \snm{Zhu}\thanksref{t3}\ead[label=e3]{htzhu@email.unc.edu}}
\and
\author{\fnms{J. S.} \snm{Marron}\thanksref{t4}
\ead[label=e4]{marron@email.unc.edu}}

\thankstext{m1}{Corresponding Author}
\thankstext{t1}{Partially supported by NSF grant DMS-0854908}
\thankstext{t2}{Partially supported by NSF grants DMS-1106912 and CMMI-0800575, NIH Challenge Grant 1 RC1 DA029425-01, and the Xerox
Foundation UAC Award}
\thankstext{t3}{Partially supported by NIH grants  RR025747-01,  P01CA142538-01,   MH086633, EB005149-01 and  AG033387}
\thankstext{t4}{Partially supported by NSF grants DMS-0606577 and DMS-0854908}
\runauthor{Dan Shen, Haipeng Shen, Hongtu Zhu and J. S. Marron}

\affiliation{University of North Carolina at Chapel Hill}

\address{Department of Biostatistics\\
University of North Carolina at Chapel Hill\\
Chapel Hill, NC 27599\\
\printead{e1}\\
\phantom{E-mail:\ }\printead*{e3}
}
\address{Department of Statistics and Operations Research\\
University of North Carolina at Chapel Hill\\
Chapel Hill, NC 27599\\
\printead{e2}\\
\phantom{E-mail:\ }\printead*{e4}
}
\end{aug}

\begin{abstract}
The aim of this paper is to establish several deep  theoretical properties of
  principal component analysis  for   multiple-component spike covariance models. Our new results reveal
    a surprising asymptotic conical structure in critical sample eigendirections under the spike models with distinguishable (or
    indistinguishable) eigenvalues,  when the sample size and/or the   number of variables (or dimension) tend to infinity.
The  consistency of the sample  eigenvectors relative to their population counterparts is determined by
 the ratio between the dimension and the product of the sample
size with the spike size.
When this ratio converges to a nonzero constant,  the sample eigenvector
converges to a cone, with  a certain angle to its corresponding  population eigenvector.
 In the High Dimension, Low Sample Size case, the angle  between the sample eigenvector and its population counterpart
converges to a  limiting distribution.
Several generalizations of the multi-spike covariance models are also explored, and additional theoretical results are presented.
\end{abstract}



\end{frontmatter}


\section{Introduction}\label{sec:introduction}

Principal Component Analysis (PCA) is  one of the most  important visualization and dimension reduction tools.
The theoretical properties of PCA,   including the sample eigenvalues, eigenvectors, and
PC scores, have been widely studied in different settings, when the sample size and/or the dimension increase to infinity. For example, Anderson (1963)~\cite{anderson1963asymptotic} studied such properties under the classical
statistical  setting with   $n\rightarrow \infty$ and
a fixed dimension $d$.
Johnstone and Lu (2009)~\cite{johnstone2009consistency} explored such properties under the random matrix setting
  with
 sample size $n\rightarrow \infty$ and  $d\sim n$.  Jung and Marron (2009)~\cite{jung2009pca}
derived such properties   in a High Dimension, Low Sample Size (HDLSS) context,
with a fixed $n$  and   $d\rightarrow\infty$. More recently, Fan et al. (2013)~\cite{fan2013large} considered scenarios where the first few leading eigenvalues increase to $\infty$ together with $d$.
See  additional theoretical results in \cite{baik2005phase,baik2006eigenvalues,onatski2006asymptotic,paul2007asymptotics,nadler2008finite,lee2010convergence,
yata2012effective,benaych2011eigenvalues,Jung2010,shen2012general} and references therein.

Generally speaking, the existing results indicate   that  the behavior of PCA strongly depend on
 the relationship among three key quantities:
the dimension, the sample size, and the spike sizes
(the relative sizes of the population eigenvalues $\{\lambda_j\}$).
For instance,   Shen et al (2012)~\cite{shen2012general} systematically investigated
the theoretical properties of  the $j$-th
sample eigenvector and eigenvalue as
 $d/(n\lambda_j)\rightarrow 0$ or $\infty$. Specifically, as $d/(n\lambda_j)\rightarrow 0$, the   $j$-th
sample eigenvector converges to the corresponding population eigenvector, whereas
 strong inconsistency follows as $d/(n\lambda_j)\rightarrow\infty$.

An interesting open  question is to investigate     the asymptotic
properties of PCA when $d/(n\lambda_j)$ converges to a constant  $c_j \in (0, \infty)$, which is the aim of this paper. A broad theoretical framework of PCA under  a broad range of cases, from the classical, through random matrix theory, and on to HDLSS, is studied here.  Firstly,  we show
a new instance of unexpected asymptotic behavior of sample eigenvectors.
 Specifically,
  the critical sample eigenvectors lie in
 a right circular cone around the corresponding population eigenvectors. Although these sample eigenvectors converge to the cone, their locations
 within the  cone are random.
 The angles of these cones have an increasing order, which is driven by  an increasing sequence of the ratios $c_j$.
   We suggest this is as surprising as the HDLSS geometric representation results discovered by Hall et al (2005)~\cite{hall2005geometric}, and further developed by Yata and Aoshima (2012)~\cite{yata2012effective}.

Secondly,  we further extend the new  results to the multi-spike cases where
the  population eigenvalues are asymptotically indistinguishable.
   We   study the angle between the corresponding sample eigenvectors
 and the subspace spanned by the indistinguishable population eigenvectors.
In HDLSS contexts,
 the cone angles are always  random variables, whereas such randomness disappears when the sample size increases.
We also show that in HDLSS settings, the PC scores are not consistent even when
the angles between the sample eigenvectors and their  population counterparts converge to 0.

Next we introduce two illustrative examples to help understand the main theoretical results in the paper,
where the eigenvalues are respectively asymptotically distinguishable (Example~\ref{example1:eigenvalues}) and indistinguishable (Example~\ref{example:angle_subspace}). Our theorems are applicable to a much broader class of general spike models.

\begin{figure}[h]
 \begin{center}
 \includegraphics[width=\textwidth]{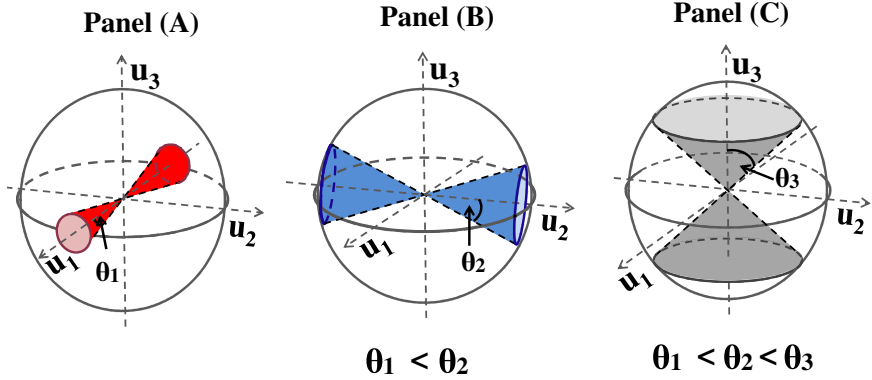}\vspace{-3.3cm}
 \end{center}
 \vspace{-.5cm}
 \caption{Geometric representation of PC directions in Example~\ref{example:angle_vector}.
 The sphere represents the space of possible sample eigenvectors. Panel (A) shows that
  the first sample eigenvector tends to  lie in the red cone, with the $\theta_1$  angle.
  Similarly, Panels (B) and (C) show
    that the second and the third sample eigenvectors respectively tend to lie in the blue and the
    gray cones, whose angles are $\theta_2$ and $\theta_3$.
    Note that the angle of the red cone is less than the blue cone, whose angle is again less than
    the gray cone.}
    \label{fig:vector}
 \end{figure}

\begin{example} (Multiple-component spike models with distinguishable eigenvalues)
Assume that $X_1,\ldots, X_n$ are random  sample vectors
from a $d$-dimensional normal distribution $N(0, \Sigma)$, where the population eigenvalues have
the following properties: as $n\rightarrow\infty$,
\begin{equation}
\left\{ \begin{array}{ll}
  \lambda_1 > \lambda_2 > \lambda_3\gg \lambda_4=\cdots=\lambda_d=1, \\
  \frac{d}{n\lambda_j}\rightarrow c_j, \quad  j=1,2,3,\; \mbox{with}\;\; 0 \leq c_1<c_2<c_3\leq\infty.
   \end{array} \right.
   \label{example1:eigenvalues}
\end{equation}

In Figure~\ref{fig:vector}, the sphere represents the space of all possible sample eigen-directions, with the first three population eigenvectors as the coordinate axes. For this particular example, our general Theorem 3.1 suggests that
\begin{itemize}

\item
As $n\rightarrow \infty$, the sample eigenvector $\hat{u}_1$ lies in the red cone, shown in Panel (A) of Fig. \ref{fig:vector},
where the angle of the cone is $\theta_1=\arccos(\frac{1}{\sqrt{1+c_1}})$. Similarly,
as $n\rightarrow\infty$, the sample eigenvectors $\hat{u}_2$ and $\hat{u}_3$
respectively lie in the blue and dark gray cones, shown in Panels (B) and (C) of Fig. \ref{fig:vector}, whereas the angles are respectively  $\theta_2=\arccos(\frac{1}{\sqrt{1+c_2}})$ and $\theta_3=\arccos(\frac{1}{\sqrt{1+c_3}})$.
Note that for $c_1<c_2<c_3$, we have $\theta_1<\theta_2<\theta_3$, as shown in Figure~\ref{fig:vector}.
\end{itemize}

In addition, our Proposition 3.1 includes the two boundary cases studied by Shen et al. (2012)~\cite{shen2012general} as special cases:
\begin{itemize}
\item
When $c_1=c_2=c_3=0$, it follows that $\theta_1=\theta_2=\theta_3=0$.
 This puts us in the domain of consistency~\cite{shen2012general}.

\item In the opposite boundary case of $c_1=c_2=c_3=\infty$, we have that $\theta_1=\theta_2=\theta_3=90$  degrees.
 This leads to  strong inconsistency~\cite{shen2012general}.
\end{itemize}

Hence, our new results  go well beyond the work of~\cite{shen2012general}, and completely characterize the transition between consistency and strong inconsistency.

\begin{figure}[h]
 \begin{center}
 \includegraphics[width=\textwidth]{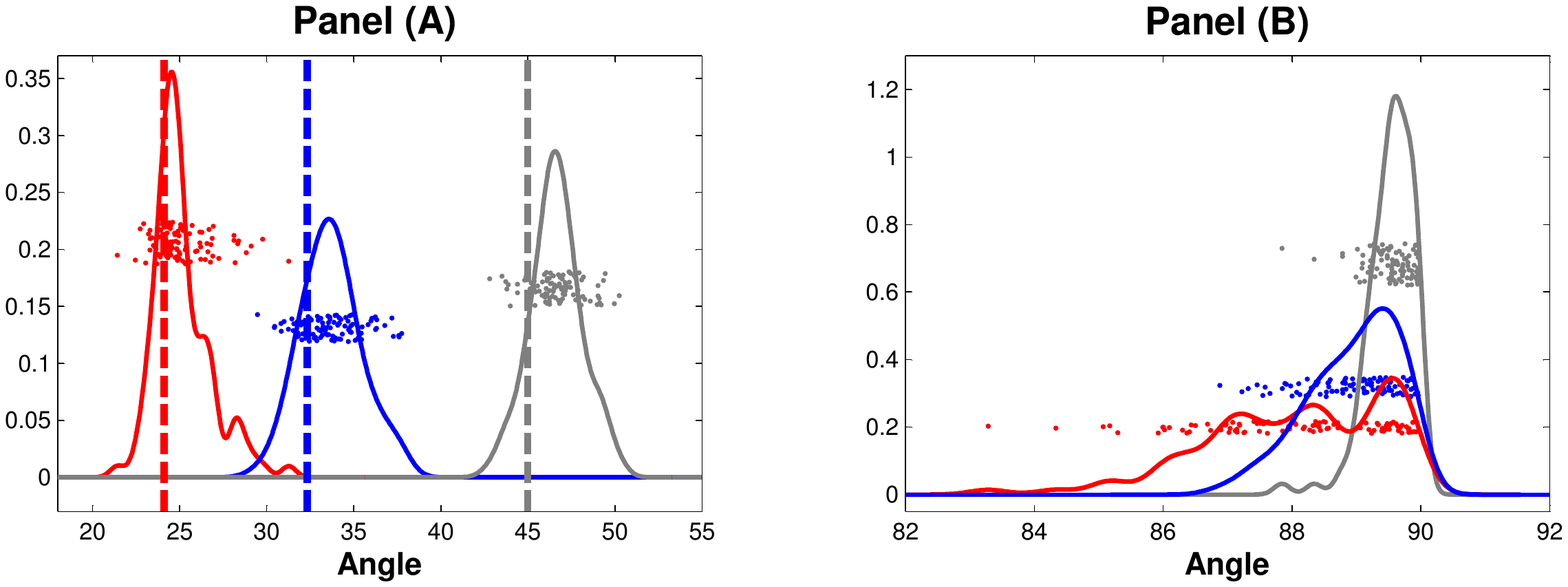}\vspace{-4cm}
 \end{center}
 \vspace{-.5cm}
 \caption{ Example~\ref{example:angle_vector}: Simulated angles between sample and population eigenvectors.
 Panel (A) shows realizations of angles between sample and population eigenvectors as colored dots (red is first, blue is second, gray is third).  Distributions are studied using kernel density estimates, and compared with the theoretical values $\theta_j$ for $j=1,2,3$, shown as dashed lines.
 Panel (B) studies randomness of eigen-directions within the cones shown in Figure~\ref{fig:vector}, by showing the distribution of pairwise angles between realizations of the sample eigenvectors.  All 3 colors are overlaid here, and all angles are very close to 90 degrees, which is very consistent with the randomness of the respective sample eigenvectors within the cones.
 }
 \label{fig:simulation}
 \end{figure}


We investigated this theoretical convergence, using simulations,
over a range of settings, with $n = 50, 100, 200, 500, 1000, 2000$, where $d/n = 50$,
and $c_1=0.2, c_2=0.4, c_3=1$. The full sequence, illustrating this convergence, is shown in Figure  
 A of the supplementary material~\cite{shen2013online}.  Figure~\ref{fig:simulation} shows the intermediate case of $n = 200$.
For one data set with this distribution,
we compute angles
between the sample and population eigenvectors. Repeating this procedure over 100 replications,
we get 100 angles for each of the first three eigenvectors, which are shown as red, blue and gray points in Panel (A).
The red, blue, gray curves are the corresponding kernel density estimates.
Panel (A) shows that the simulated angles are very close to the corresponding theoretical angles
$\theta_j$, $j=1,2,3$, shown as dashed vertical lines.

Panel (B) in Figure~\ref{fig:simulation} studies randomness of eigen-directions within the cones
shown in Figure~\ref{fig:vector}.
We calculate pairwise angles between realizations of the sample eigenvectors for the three cones,
showing angles and kernel density estimates using colors as in Panel (A) of Figure~\ref{fig:simulation}.
All angles are very close to 90 degrees, which is  consistent
with randomness in high dimensions, see~\cite{hall2005geometric, yata2012effective, jung2009pca,Jung2010} and the more recent work of Cai et al. (2013)~\cite{cai2013distributions}.
In fact, the regions represented by circles in Figure 1, are actually $d-1$ dimensional hyperspheres, so
 the sample eigenvectors should be thought of as d-1 dimensional as $d, n\rightarrow \infty$.

\label{example:angle_vector}
\end{example}

\begin{example} (Multiple-component spike models with indistinguishable eigenvalues)
We again assume that $X_1,\ldots, X_n$ are random sample vectors
from a $d$-dimensional normal distribution $N(0, \Sigma)$. Different from Example~\ref{example:angle_vector}, the six leading
population eigenvalues of $\Sigma$ fall into three asymptotically separable pairs as follows: as $n\rightarrow \infty$
\begin{equation*}
\left\{ \begin{array}{ll}
  \lambda_1=\lambda_2 > \lambda_3=\lambda_4 > \lambda_5=\lambda_6\gg \lambda_7=\cdots=\lambda_d=1, \\
  \frac{d}{n\lambda_{2j-1}}\rightarrow c_j, \quad  j=1,2,3,\; \mbox{with}\;\; 0 \leq c_1<c_2<c_3\leq\infty.
   \end{array} \right.
\end{equation*}

\begin{figure}[h]
 \begin{center}
 \includegraphics[width=\textwidth]{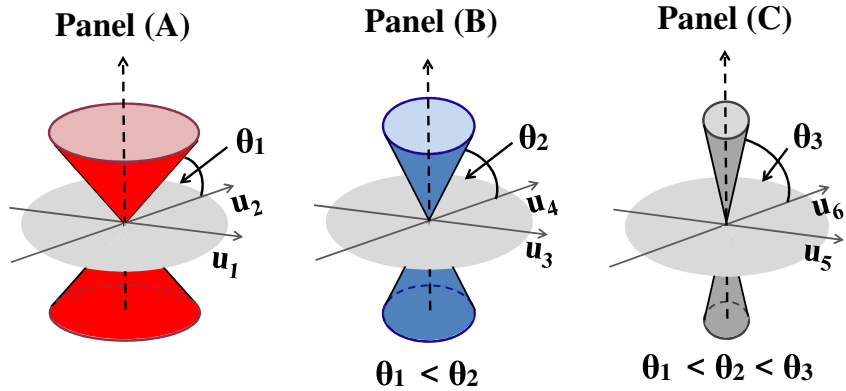}\vspace{-3.3cm}
 \end{center}
 \vspace{-.5cm}
 \caption{Example~\ref{example:angle_subspace}: Geometric representation of PC directions.
 Panel (A) shows the cone to which
   the first group of  sample eigenvectors converge in the red. This cone has angle $\theta_1$  with the
 gray subspace, generated by the first group of population eigenvectors.
  Similarly, Panel (B) (Panel (C)) shows  the cone to which the second (third) group of
   sample eigenvectors converges  shown as a blue (dark gray) cone, which has angle  $\theta_2$ ($\theta_3$) with the subspace, generated by the second (third) group of  population eigenvectors.}
 \label{fig:subspace}
\end{figure}

Our general Theorem 3.2, when applied to the current example, reveals the following insights:
\begin{itemize}
\item  Panel (A) in Figure~\ref{fig:subspace} shows, as
a red cone, the region where the first group of sample eigenvectors $\hat{u}_1$ and $\hat{u}_2$  lie in the limit
as $n\rightarrow\infty$.
This has the angle $\theta_1=\arccos(\frac{1}{\sqrt{1+c_1}})$  with the gray subspace,
generated by the first group of population eigenvectors $u_1$ and $u_2$.
Similarly, Panel (B) (Panel (C)) presents, as a blue (gray) cone, the
region where the second (third) group of sample eigenvectors $\hat{u}_3$ and $\hat{u}_4$ ($\hat{u}_5$ and $\hat{u}_6$) lie
 in the limit as $n\rightarrow \infty$. This has the angle $\theta_2=\arccos(\frac{1}{\sqrt{1+c_2}})$ ($\theta_3=\arccos(\frac{1}{\sqrt{1+c_3}})$)  with
 the  subspace, generated by  the second (third) group of population eigenvectors
 $u_3$ and $u_4$ ($u_5$ and $u_6$). Note that for $c_1<c_2<c_3$, we have $\theta_1<\theta_2<\theta_3$, as shown in Figure~\ref{fig:subspace}.
\end{itemize}

Furthermore, our Proposition B.1 in the supplementary document~\cite{shen2013online} considers boundary cases of our general framework, which includes the results of Shen et al. (2012)~\cite{shen2012general} as special cases:
\begin{itemize}
\item For $c_1=c_2=c_3=0$, it follows that $\theta_1=\theta_2=\theta_3=0$.
    This puts us in the domain of subspace consistency, as studied in Theorem 4.3 of~\cite{shen2012general}.

\item When $c_1=c_2=c_3=\infty$, we have that $\theta_1=\theta_2=\theta_3=90$ degrees. This leads to strong inconsistency, as studied in Theorem 4.3 of~\cite{shen2012general}.
\end{itemize}
\label{example:angle_subspace}
\end{example}

The rest of the paper is organized as follows.
Section~\ref{sec:assumption} introduces the assumptions and notation
relevant to the theorems in the paper.
 Section~\ref{sec:hdmss} studies
the asymptotic properties of PCA
for  multiple spike models  with  distinguishable (or indistinguishable) eigenvalues
as $n\rightarrow\infty$.
Section~\ref{sec:hdlss} studies
the asymptotic properties of PCA in the HDLSS contexts.
Section~\ref{proofs} contains the technical proofs of the main theorems.
Additional simulation studies and proofs can be found in the supplementary document~\cite{shen2013online}.

\section{ Assumptions and Notation}\label{sec:assumption}
Let $X_1,\ldots, X_n$ be random vectors from
a $d$-dimensional normal distribution $N(\xi, \Sigma)$, where $\xi$ is a $d\times 1$ mean vector and $\Sigma$ is a $d\times d$ covariance matrix.
Let $\{(\lambda_k, u_k): k=1, \cdots, d\}$ be the eigenvalue-eigenvector pairs of $\Sigma$ such that $\lambda_1\geq\lambda_2\geq\ldots\geq\lambda_d>0$. Thus,
  $\Sigma$ has the following eigen-decomposition
$$
\Sigma=U\Lambda U^T,
$$ where $\Lambda=\mbox{diag}(\lambda_1,  \ldots, \lambda_d)$ and
$U=[u_1,\ldots,u_d]$.  Since the relative sizes, rather than the absolute values, of the population eigenvalues
  affect the asymptotic properties of PCA,   we assume that $\lambda_d=1$ throughout the rest of the paper.


Let $\overline{X}$ be the sample mean.  As discussed in~\cite{paul2007augmented},
\begin{equation}
\sum_{i=1}^n (X_i-\overline{X})(X_i-\overline{X})^T\quad \mbox{has\;the\;same\;distribution\;as}\quad \sum_{i=1}^{n-1} Y_iY_i^T,
\label{centered-CovarinaceMatrix}
\end{equation}
where $Y_i$ are i.i.d  random vectors from $N(0,\Sigma)$.
It follows from~\eqref{centered-CovarinaceMatrix} that the sample covariance matrix is location invariant.
 Thus,  we can assume without loss of generality (WLOG):
\begin{assumption}
$X_1,\ldots, X_n$ are i.i.d
random vectors from a $d$-dimensional normal distribution $N(0,
\Sigma)$.
\label{gauss-assumption}
\end{assumption}
Denote the $j$th normalized population PC score vector as
\begin{equation}
S_j=(S_{1,j},\cdots, S_{n,j})^T=\lambda_j^{-\frac{1}{2}}(u_j^T X_1,\cdots,  u_j^TX_n)^T, \quad j=1,\cdots, d,
\label{PCs_population}
\end{equation}
and define $Z$ as the $n\times d$ random matrix as
\begin{equation}
Z=(z_{i,j})_{n\times d}= X^T U \Lambda^{-\frac{1}{2}},
\label{def:z}
\end{equation}
where $X=[X_1,\ldots,X_n]$  and $z_{i,j}$, $i=1,\cdots,n$, $j=1,\cdots,d$ are i.i.d random variables from $N(0,1)$.

Let $\{(\hat \lambda_k, \hat u_k): k=1, \cdots, d\}$ be the eigenvalue-eigenvector pairs of the sample covariance matrix $\hat{\Sigma}=n^{-1}XX^{T}$ such that $\hat \lambda_1\geq\hat \lambda_2\geq\ldots\geq \hat \lambda_d$. Thus,
 $\hat{\Sigma} $ can  be
decomposed as
\begin{equation}
\hat{\Sigma}=\hat{U}\hat{\Lambda}\hat{U}^T,
\label{sample:covariance}
\end{equation}
where  $\hat{\Lambda}=\mbox{diag}(\hat{\lambda}_1,\ldots,\hat{\lambda}_d)$
and $\hat{U}=[\hat{u}_1,\ldots,\hat{u}_d]$.
Note that the data matrix $n^{-\frac{1}{2}}X$ has the singular value decomposition such that  $n^{-\frac{1}{2}}X=\sum_{j=1}^d \hat{\lambda}_j^{\frac{1}{2}}\hat{u}_j\hat{v}^T_j$, where $\hat{v}_j=(\hat{v}_{1,j},\cdots,\hat{v}_{n,j})^T$ for  $j=1, \cdots, d$. Thus,  the $j$th normalized sample PC score vector is  given by 
\begin{equation}
\hat{S}_j=(\hat{S}_{1,j},\cdots, \hat{S}_{n,j})^T=(\hat{v}_{1,j},\cdots,\hat{v}_{n,j})^T, \quad j=1,\cdots, d.
\label{PCs_sample}
\end{equation}

{{We introduce an asymptotic notation.
Assume that $\{\xi_{k}: k=1,\ldots,\infty\}$ ($k=n$ or $d$)
 is a sequence of random variables  and
$\{e_{k}: k=1,\ldots,\infty\}$ is a sequence of constants.
 Denote $\xi_{k}={\rm O}_{\rm a.s}\left(e_{k}\right)$ if $\overline{\mbox{lim}}_{k\rightarrow \infty}\left|\frac{\xi_{k}}{e_{k}}\right|\leq \zeta$ almost surely with $P(0< \zeta < \infty)=1$ .
}}

\section{Growing sample size asymptotics}
\label{sec:hdmss}
We now
 study asymptotic properties of PCA as $n\rightarrow\infty$.
 We consider
 multiple component spike models with distinguishable
 population eigenvalues  in Section~\ref{subsec:distinct-eigenvalues-hdmss}  and
  with
indistinguishable eigenvalues in Section~\ref{subsec:msm-tiered-eigenvalues}.
Moreover, we vary   $d$   from the classical $d$ fixed asymptotics, through the random matrix version with $d \sim n$,
 all the way to the high dimension medium sample size (HDMSS) asymptotics of Cabanski et al (2010)~\cite{cabanski2010swiss} and Yata and Aoshima (2012)~\cite{yata2012inference} with
 $d\gg n\rightarrow\infty$.

\subsection{Multiple component spike models with distinguishable eigenvalues}
\label{subsec:distinct-eigenvalues-hdmss}

We consider  multiple component spike models with $m$
dominating spikes where finite $m\in [1,n\wedge d]$.
The population eigenvalues are assumed to satisfy the following two assumptions:
\begin{enumerate}
\item[$\mathcal{A}1$.] As $n \rightarrow \infty$, 
$\lambda_1 >\cdots > \lambda_m\gg \lambda_{m+1}\rightarrow\cdots \rightarrow\lambda_d=1.$

\item[$\mathcal{A}2$.]
$
{\rm As}\; n\rightarrow \infty,\quad\frac{d}{n\lambda_j}\rightarrow c_j,\quad {\rm where}\quad 0< c_1<\cdots< c_m <\infty.
$

\end{enumerate}

We first make several comments about Assumptions $\mathcal{A}1$ and $\mathcal{A}2$.

\begin{itemize}
\item Assumption $\mathcal{A}1$ includes two separate parts:

\begin{enumerate}

\item[(a)] The $\lambda_1 >\cdots > \lambda_m$  part makes it possible to separately consider the first $m$ principle component signals and study the corresponding asymptotic properties.

\item[(b)] The $\lambda_m\gg \lambda_{m+1}\rightarrow\cdots \rightarrow\lambda_d=1$ enables clear separation of the signal (contained in the first $m$ components) from the noise (in the higher order components),
  which then helps to derive the asymptotic properties of the first $m$ sample eigenvalues, eigenvectors, and PC scores.

\end{enumerate}

\item Assumption  $\mathcal{A}2$  is  the critical case, in which  the  \emph{positive} information and the \emph{negative}  are of the same order. In particular,  increasing $n$ and the spike  positively impacts the consistency of  PCA,  whereas increasing $d$  has a negative impact.
\end{itemize}

\begin{figure}[h]
 \begin{center}
 \includegraphics[width=\textwidth]{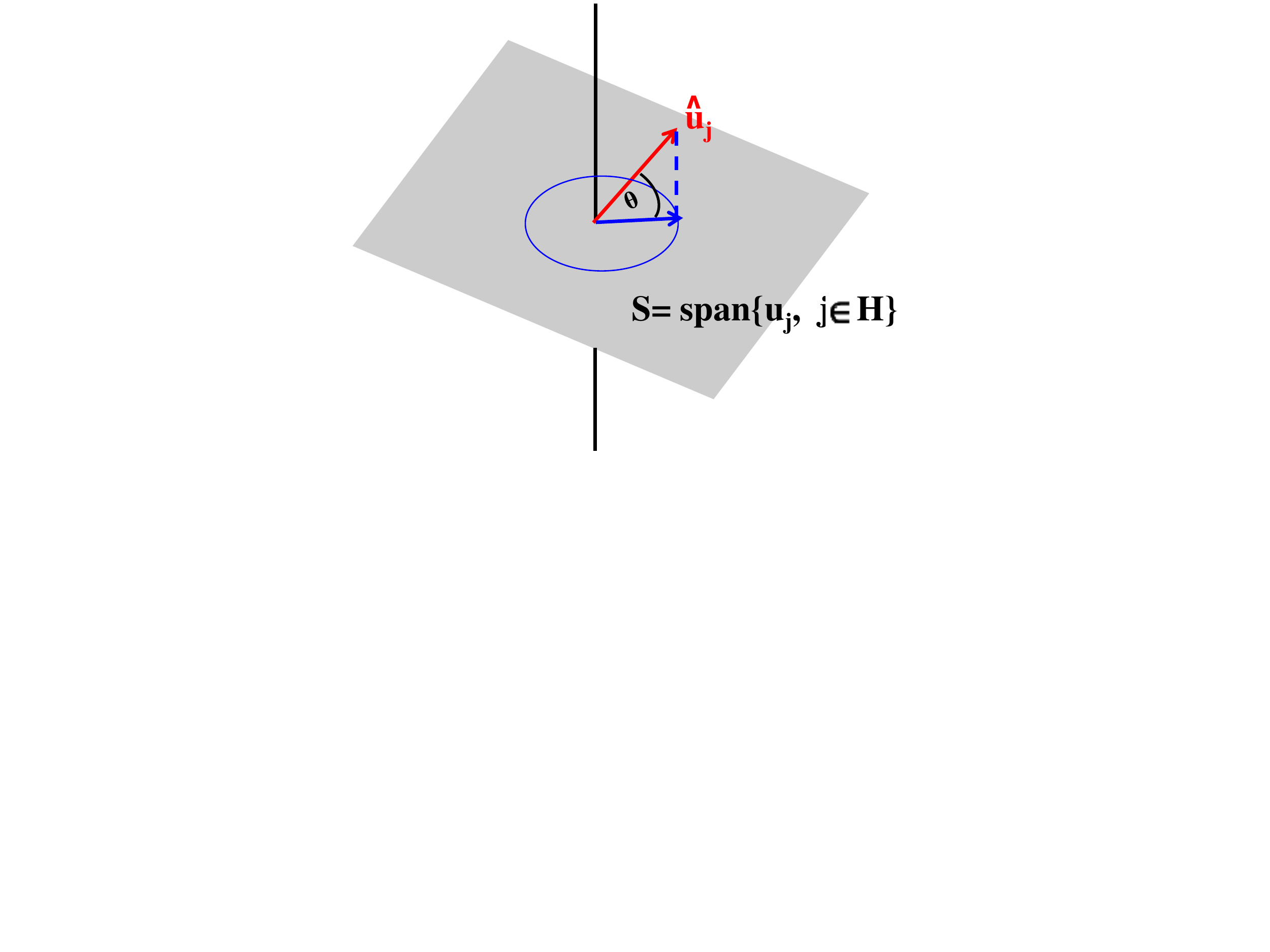}\vspace{-4.3cm}
 \end{center}
 \vspace{-.5cm}
 \caption{Angle between the sample eigenvector $\hat{u}_j$ and the space $\mathbb{S}$. The blue vector is
 the projection of the red vector $\hat{u}_j$ onto the space $\mathbb{S}$.}
 \label{fig:angle_definition}
\end{figure}

While the main focus of our results is the signal eigenvectors, some notation for the noise eigenvectors is also useful.
According to Assumption $\mathcal{A}1$, the noise sample eigenvalues
whose indices are greater than $m$ can not be asymptotically distinguished, so the
corresponding eigenvectors should be treated as a whole.
Therefore, we define the noise index set $H=\{m+1,\cdots,d\}$, and denote the space spanned by these noise eigenvectors as
\begin{equation}
\mathbb{S}=\mbox{span}\{u_j, j \in H\}.
\label{subspace:SS}
\end{equation}
For each sample eigenvector $\hat{u}_j$, $j\in H$, we study the angle between $\hat{u}_j$ and the space $\mathbb{S}$,
as defined in~\cite{jung2009pca,shen2012general} and illustrated in Figure~\ref{fig:angle_definition}, i.e.
the angle between $\hat{u}_j$ (the red vector) and its projection onto $\mathbb{S}$ (the blue vector).


The following theorem derives the asymptotic properties of the first
$m$ sample eigenvalues and eigenvectors. In addition, the theorem also
 shows that, for $j=m+1,\cdots,[n\wedge d]$,
the angle between $\hat{u}_j$ and $u_j$
 goes to 90 degrees, whereas the angle between $\hat{u}_j$ and the space
$\mathbb{S}$ goes  to 0.

\begin{theorem}  Under Assumptions~\ref{gauss-assumption}, $\mathcal{A}1$,
and $\mathcal{A}2$,
 as $n \rightarrow \infty$, the sample eigenvalues satisfy
\begin{equation}
\left\{ \begin{array}{ll}
  \frac{\hat{\lambda}_j}{\lambda_j} \xrightarrow{\rm a.s} 1+c_j,\; 1\leq j\leq m, \\
  \frac{n\hat{\lambda}_j}{d\lambda_j} \xrightarrow{\rm a.s} 1, \; m+1\leq j\leq [n\wedge d],\\
   \end{array} \right.
\label{eigenvalue:distinct}
\end{equation}
and the sample eigenvectors satisfy
\begin{equation}
\left\{\begin{array}{lll}
\mid<\hat{u}_j, u_j>\mid \xrightarrow{a.s} (1+c_j)^{-\frac{1}{2}},\quad 1\leq j\leq m,\\
\mid<\hat{u}_j, u_j>\mid={\rm O}_{\rm a.s} \left\{(\frac{n}{d})^{\frac{1}{2}}\right\}, \quad  m+1\leq j\leq [n\wedge d],\\
\mbox{angle}<\hat{u}_j, \mathbb{S}>\xrightarrow{a.s} 0, \quad m+1\leq j\leq [n\wedge d].
\end{array} \right.
\label{eigenvector:distinct}
\end{equation}
\label{Th:distinct}
\end{theorem}

We now offer several remarks regarding Theorem~\ref{Th:distinct}.

\begin{remark}
The results of~\eqref{eigenvalue:distinct} and~\eqref{eigenvector:distinct} suggest that, as the  eigenvalue index increases, the proportional bias
between the sample and population eigenvalue increases, so the angle between the sample and corresponding population eigenvectors increases.
 This is because larger eigenvalues (i.e. with small indices) contain more positive information, which makes the corresponding sample eigenvalues/eigenvectors less biased. These results are graphically illustrated in Figure 1 and empirically verified in Figure 2, for the specific model in Example 1.1. More empirical support is provided in the supplementary material~\cite{shen2013online}.
\end{remark}

\begin{remark} Theorem~\ref{Th:distinct} can be extended to include the classical and random matrix cases, by allowing $c_{m_0}=0$ for some $m_0\leq m$, which suggests that positive information dominates in the leading $m_0$ spikes.
Then Assumptions $\mathcal{A}1$ and $\mathcal{A}2$ respectively become
\begin{enumerate}
\item[$\mathcal{A}3$.] as $n \rightarrow \infty$, the population eigenvalues satisfy
\begin{equation*}
\lambda_1 >\cdots > \lambda_{m_0}\gg \lambda_{m_0+1}>\cdots> \lambda_{m} \gg \lambda_{m+1}\rightarrow\cdots \rightarrow\lambda_d=1.
\end{equation*}

\item[$\mathcal{A}4$.] as $n\rightarrow \infty$, $d/(n\lambda_j)\rightarrow c_j$ for $j=1,\cdots,m$,
where $0=c_1=\cdots= c_{m_0}< c_{m_0+1}<\cdots< c_m <\infty$.

\end{enumerate}

For the classical case with fixed dimension $d$, $m_0=m=d$ in Assumptions $\mathcal{A}3$ and $\mathcal{A}4$.
For random matrix cases with $n\sim d$, $m_0=m$ in Assumptions $\mathcal{A}3$ and $\mathcal{A}4$.
 Since $c_1=\cdots= c_{m_0}=0$ in Assumption $\mathcal{A}3$, if the eigenvalue index is less than or equal to $m_0$, the corresponding sample eigenvalues and eigenvectors are consistent.
 These results are summarized in the following Proposition~\ref{corr:distinct}(a).
 \end{remark}

 \begin{remark}
Another extension of Theorem~\ref{Th:distinct} is to allow $c_{m_0+1}=\infty$ for some $m_0\le m$, i.e. negative information dominates in higher-order spikes. This contains the HDMSS cases~\citep{cabanski2010swiss, yata2012inference}, where $d\gg n\rightarrow\infty$.
Assumption $\mathcal{A}1$ then becomes Assumption $\mathcal{A}3$, and Assumption $\mathcal{A}2$  becomes
\begin{enumerate}
\item[$\mathcal{A}5$.] as $n\rightarrow \infty$, $d/(n\lambda_j)\rightarrow c_j$ for $j=1,\cdots,m$,
where $0< c_1<\cdots< c_{m_0}< c_{m_0+1}=\cdots= c_m =\infty$.

\end{enumerate}

Since $c_{m_0+1}=\cdots= c_m =\infty$, for index $j\ge m_0+1$,  the proportional error between the sample and population
eigenvalues goes to infinity, and the angle between the corresponding  sample and population eigenvectors converges to 90 degrees.
 These results are summarized in Proposition~\ref{corr:distinct}(b).
\end{remark}

\begin{proposition}


 \begin{enumerate}

\item[(a)] Under Assumptions~\ref{gauss-assumption}, $\mathcal{A}3$ and $\mathcal{A}4$, the sample eigenvalues and eigenvectors satisfy
\begin{equation*}
  \hat{\lambda}_j/\lambda_j \quad \mbox{and} \quad
  \mid <\hat{u}_j, u_j> \mid  \quad \xrightarrow{a.s} 1, \quad 1\leq j\leq m_0,
  \label{proposition:case1}
  \end{equation*}
     and the properties of the other sample eigenvalues  and eigenvectors remain the same
    as in Theorem~\ref{Th:distinct}.

\item[(b)] Let $H=\{m_0+1,\cdots,d\}$ and define
  $\mathbb{S}$ as in~\eqref{subspace:SS}. If Assumption $\mathcal{A}4$ in (a) is
  replaced by Assumption $\mathcal{A}5$,
  the sample eigenvalues satisfy
  \begin{equation*}
   n\hat{\lambda}_j/d \xrightarrow{\rm a.s} 1, \quad m_0+1\leq j\leq m,
   \label{proposition:case1:eigenvalues}
   \end{equation*}
   and the sample eigenvectors satisfy
\begin{equation*}
  \left\{ \begin{array}{ll}
   \mid<\hat{u}_j, u_j>\mid ={\rm O}_{\rm a.s}\left\{(\frac{n\lambda_j}{d})^{\frac{1}{2}}\right\},\\
  \mbox{angle}<\hat{u}_j, \mathbb{S}>\xrightarrow{a.s} 0,\\
  \end{array} \right.
   \quad m_0 \leq j\leq [n\wedge d];
   \label{proposition:case2:eigenvectors}
   \end{equation*}
the properties of the other sample eigenvalues  and eigenvectors remain the same as in Theorem~\ref{Th:distinct}.

\item[(c)] In addition, if Assumption $\mathcal{A}4$ in (a) is strengthened to $d/\lambda_{m_0}\rightarrow 0$,
then the sample PC scores satisfy
\begin{equation*}
\left|\frac{\hat{S}_{i,j}}{S_{i,j}}\right|\xrightarrow{a.s} 1, \quad i=1,\cdots n, \;
j=1,\cdots, m_0.
\label{PCs:distinct}
\end{equation*}

\end{enumerate}

\label{corr:distinct}

\end{proposition}

%
%

\subsection{Multiple component spike models with indistinguishable eigenvalues}
\label{subsec:msm-tiered-eigenvalues}
We now consider spike models with the $m$ leading eigenvalues being  grouped into $r$($\ge1$) \emph{tiers},
each of which contains eigenvalues that are either the same or have the same limit.
The eigenvalues within different tiers have different limits.
Specifically, the first $m$ eigenvalues are grouped into $r$ tiers, in which there are $q_k$ eigenvalues in the $k$th tier such that $\sum_{l=1}^r q_l=m$. Define $q_0=0$, $q_{r+1}=d-\sum_{l=1}^r q_l$, and the index set of the eigenvalues in the $k$th tier as
\begin{equation}
H_k=\left\{\sum_{l=0}^{k-1} q_l+1,\sum_{l=0}^{k-1}
q_l+2,\cdots, \sum_{l=0}^{k-1} q_l+q_k\right\}, \quad k=1,\cdots,r+1.
\label{index}
\end{equation}

We make the following assumptions on the tiered eigenvalues:
\begin{enumerate}
\item[$\mathcal{B}1$.] The eigenvalues  in the $k$th tier
have the same limit $\delta_k(>0)$:
\begin{equation*}
\mbox{lim}_{n\rightarrow \infty} \frac{\lambda_j}{\delta_k}=1, \quad j\in H_k, \;k=1,\cdots,r.
\end{equation*}

\item[$\mathcal{B}2$.] The eigenvalues in different tiers
have different limits:
\begin{equation*}
{\rm as}\; n \rightarrow \infty,\; \quad \delta_1 > \cdots > \delta_r\gg \lambda_{m+1}\rightarrow\cdots\rightarrow\lambda_d=1.
\end{equation*}

\item[$\mathcal{B}3$.] The ratio between the dimension and the product of the sample size with
  eigenvalues in the same tier converges to a constant:
\begin{equation*}
{\rm as}\; n\rightarrow \infty, \quad
\frac{d}{n\delta_k}\rightarrow c_k,\quad {\rm with} \quad 0<c_1<\cdots<c_r<\infty.
\end{equation*}
\end{enumerate}

  Assumptions  $\mathcal{B}2$ and $\mathcal{B}3$ are natural extensions of Assumptions $\mathcal{A}1$ and $\mathcal{A}2$. In Assumption  $\mathcal{B}2$,   the signal contained in the first $r$ tiers of eigenvalues is well separated from the noise,  and hence the asymptotic properties of the sample eigenvalues and eigenvectors in the first $r$ tiers can be obtained. Assumption $\mathcal{B}3$ suggests that
 the positive information (sample size and spike size) and the negative information (dimension)
 are of the same order.

Since the sample eigenvalues within the same tier can not be asymptotically identified,
 the corresponding sample eigenvectors are indistinguishable.
For $j\in H_k$, in order
to study the asymptotic properties of the sample eigenvector $\hat{u}_j$,
we consider the angle between $\hat{u}_j$ and the subspace spanned by the population eigenvectors $u_j$ in the same tier, defined as
 \begin{equation}
 \mathbb{S}_k=\mbox{span}\{u_j, j \in H_k\}.
 \label{subspace:S}
 \end{equation}

Our theoretical results are summarized in the following theorem.

\begin{theorem}  Under Assumptions~\ref{gauss-assumption},
$\mathcal{B}1$, $\mathcal{B}2$ and $\mathcal{B}3$, as
$n \rightarrow \infty$,
the sample eigenvalues satisfy
\begin{equation}
\left\{ \begin{array}{ll}
  \frac{\hat{\lambda}_j}{\lambda_j} \xrightarrow{\rm a.s} 1+c_k,\quad j\in H_k, \; k=1,\cdots, r, \\
  \frac{n\hat{\lambda}_j}{d\lambda_j} \xrightarrow{\rm a.s} 1, \quad m+1\leq j\leq [n\wedge d],\\
   \end{array} \right.
\label{eigenvalue:indistinct}
\end{equation}
and the sample eigenvectors satisfy
\begin{equation}
\left\{\begin{array}{lll}
\mbox{angle}<\hat{u}_j, \mathbb{S}_k>
\xrightarrow{a.s} \arccos\left\{(1+c_k)^{-\frac{1}{2}}\right\},\;j\in H_k, \; k=1,\cdots,r,\\
\mid<\hat{u}_j, u_j>\mid={\rm O}_{\rm a.s}\left\{(\frac{n}{d})^{\frac{1}{2}}\right\}, \quad  m+1\leq j\leq [n\wedge d],\\
\mbox{angle}<\hat{u}_j, \mathbb{S}_{r+1}>
\xrightarrow{a.s} 1, \quad m+1\leq j\leq [n\wedge d].
\end{array} \right.
\label{eigenvector:indistinct}
\end{equation}

\label{Th:multiple-tiered}
\end{theorem}

 Theorem~\ref{Th:multiple-tiered} is an extension of
  Theorem~\ref{Th:distinct}. For higher-order
eigenvalues, the sample eigenvalues are  more biased, while  the angles between the sample eigenvectors and the subspaces spanned by their population counterparts in the same tiers are larger. See Figure 3 for an illustration of   the specific model considered in Example 1.2.  Theorem~\ref{Th:multiple-tiered} can be extended to cover the classical, random matrix, and HDMSS cases,
which is done in Section B of the supplementary material~\cite{shen2013online}.

\section{High dimension, low sample size asymptotics}
\label{sec:hdlss}

We now  study the asymptotic properties of PCA in the HDLSS context.  In this case, the ratios between the sample eigenvalues and their population counterparts converge to non degenerate random variables, as do
the angles between the sample eigenvectors and the space spanned by the corresponding population eigenvectors. This phenomenon of random limits does not exist when $n$ increases to $\infty$  as shown in Section~\ref{sec:hdmss}.

Since the sample size   is fixed, we can not distinguish  the two types of spike models considered respectively in Sections 3.1 and 3.2. Hence, we merge the model assumptions there into the following corresponding assumptions:
\begin{enumerate}
\item[$\mathcal{C}1$.] For fixed $n$, as $d \rightarrow \infty$, $
\lambda_1 \geq\cdots \geq \lambda_m\gg \lambda_{m+1}\rightarrow\cdots \rightarrow\lambda_d=1$.

\item[$\mathcal{C}2$.] For fixed $n$, as $d\rightarrow \infty$,
\begin{equation*}
\frac{d}{n\lambda_j}\rightarrow c_j,\quad\quad
{\rm with}\quad 0< c_1\leq\cdots\leq c_m <\infty.
\end{equation*}
\end{enumerate}
In particular, Assumption $\mathcal{C}1$  is parallel to Assumptions $\mathcal{A}1$, $\mathcal{B}1$ and
$\mathcal{B}2$, while Assumption $\mathcal{C}2$ corresponds to Assumptions $\mathcal{A}2$ and $\mathcal{B}3$.

As stated below in Theorem~\ref{Th:multiple-tiered:hdlss}, the sample eigenvalues
and eigenvectors converge to non-degenerate random variables rather than constants. We   define several quantities in order to describe the limiting random variables.
Define the $m\times d$  matrix
 \begin{equation*}
\mathbb{M}=[\mathbb{C},  0_{m \times (d-m)}]_{m\times d},
 \quad
 \label{CC:matrix}
\end{equation*}
where $\mathbb{C}=\mbox{diag}\{{c^{-1/2}_1},\cdots, {c^{-1/2}_m}\}$
 is an $m\times m$ diagonal matrix and $0_{m \times (d-m)}$ is the $m\times (d-m)$ zero matrix.
In addition, define the random matrix $\mathcal{W}$ as
\begin{equation}
\mathcal{W}=\mathbb{M} Z^TZ \mathbb{M}^T,
\label{wishart}
\end{equation}
where $Z$ is defined in (2.3).
The eigenvalues of the random matrix $\mathcal{W}$ appear in the random limits of Theorem~\ref{Th:multiple-tiered:hdlss}, as in~\eqref{eigenvalue:indistinct:hdlss} and~\eqref{eigenvector:indistinct:hdlss}.

Given the fixed sample size, the sample eigenvalues can not be asymptotically distinguished,
 nor can the corresponding sample eigenvectors. To study the asymptotic behavior of the sample eigenvectors, we need to consider the space $\mathbb{S}_k$ spanned by the corresponding population eigenvectors, as defined in~\eqref{subspace:S}, with the two index sets being $H_1=\left\{1,\cdots,m\right\}$ and $H_2=\left\{m+1,\cdots,d\right\}$.


We are now ready to state the main theorem in the HDLSS contexts.

\begin{theorem}  Under Assumptions 2.1,
$\mathcal{C}1$ and $\mathcal{C}2$,  for fixed $n$, as
$d \rightarrow \infty$,
the sample eigenvalues satisfy
\begin{equation}
\left\{ \begin{array}{ll}
  \frac{\hat{\lambda}_j}{\lambda_j} \xrightarrow{\rm a.s} \frac{c_j}{n}\lambda_j(\mathcal{W})+c_j,\quad 1\leq j\leq m, \\
  \frac{n\hat{\lambda}_j}{d\lambda_j} \xrightarrow{\rm a.s} 1, \quad m+1\leq j\leq n,\\
   \end{array} \right.
\label{eigenvalue:indistinct:hdlss}
\end{equation}
where $\mathcal{W}$ is defined in~\eqref{wishart}, and the sample eigenvectors satisfy
\begin{equation}
\left\{\begin{array}{lll}
\mbox{angle}<\hat{u}_j, \mathbb{S}_1>
\xrightarrow{a.s} \arccos\left\{ \left( 1+\frac{n}{\lambda_j(\mathcal{W})}\right)^{-\frac{1}{2}}\right\},\quad 1\leq j \leq m,\\
\mid<\hat{u}_j, u_j>\mid={\rm O}_{\rm a.s}(d^{-\frac{1}{2}}), \quad  m+1\leq j\leq n,\\
\mbox{angle}<\hat{u}_j, \mathbb{S}_2>
\xrightarrow{a.s} 1, \quad m+1\leq j\leq n.
\end{array} \right.
\label{eigenvector:indistinct:hdlss}
\end{equation}

\label{Th:multiple-tiered:hdlss}
\end{theorem}

Three remarks are offered below regarding Theorem~\ref{Th:multiple-tiered:hdlss}.

\begin{remark}
If $m=1$ in Theorem~\ref{Th:multiple-tiered:hdlss}, i.e. for single-component spike models, then the first sample eigenvalue and eigenvector
satisfy
\begin{equation*}
  \left\{ \begin{array}{ll}
  \frac{\hat{\lambda}_1}{\lambda_1} \xrightarrow{a.s} \frac{\chi_n^2}{n}+c_1,\\
  \mid <\hat{u}_1, u_1>\mid \xrightarrow{a.s} \left( 1+\frac{n c_1}{\chi_n^2}\right)^{-\frac{1}{2}},
  \end{array} \right.
\end{equation*}
where $\chi_n^2$ is the Chi-square distribution with $n$ degrees of freedom.
This result is consistent with Theorem 1 of Jung et al. (2012)~\cite{Jung2010}.
\end{remark}

\begin{remark}
For $1\leq j \leq m$, as the relative size of the eigenvalue decreases,
the angle between $\hat{u}_j$ and $\mathbb{S}_1$  increases. However, this phenomenon is not as strong as in the growing sample size settings studied in Section 3, where
the sample eigenvectors can be separately studied, and the corresponding angles have a non-random increasing order.
\end{remark}

\begin{remark}
Assumption $\mathcal{C}2$ can be relaxed to include boundary cases, in which there exists an integer $m_0\in[1, m]$ such that $c_{m_0}=0$, i.e. positive information dominates in the leading $m_0$ spikes; or $c_{m_0+1}=\infty$, i.e. negative information dominates in the remaining high-order spikes. These theoretical results are presented in Section C of the
supplementary material~\cite{shen2013online}.
\end{remark}

\section{Proofs}\label{proofs}

We now provide some proofs of our theorems
as $n\rightarrow\infty$.
For the sake of space, we only present detailed proof for the properties of the sample
eigenvectors here, which is the most challenging part.
In contrast to showing consistency or inconsistency of the sample eigenvector,
this proof requires precise calculation of the degree
of inconsistency, i.e. the limiting angles
between the sample and population eigenvectors.
We relegate the derivations regarding the sample eigenvalues to Section D of
the supplementary material~\cite{shen2013online}, which also contains proofs of
Proposition 3.1, Theorem 4.1, as well as extensions of Theorems 3.2 and 4.1.

The critical ideas of the proof are to first partition the sample eigenvector matrix $\hat{U}$ into
 sub-matrices, corresponding to the group index $H_k$. Then through careful analysis,
we explore the connections between sample eigenvectors and eigenvalues and then use the sample eigenvalue properties
 to study the asymptotic properties of the sample eigenvectors.

WLOG, we  assume that $\lambda_{m+1}=\cdots=\lambda_d=1$.
Due to the invariance property of the angle between the sample and population eigenvectors,
see Shen et al. (2012)~\cite{shen2012general},
 we assume WLOG that the population eigenvectors $u_j=e_j$, $j=1,\ldots, d$,
where the $j$-th component of $e_j$ equals 1 and the rest are zero.
It follows that the inner product between the sample and population eigenvectors
satisfies
\begin{equation*}
\mid<\hat{u}_j,u_j>\mid^2=\mid<\hat{u}_j,e_j>\mid^2=\hat{u}^2_{j,j},
\label{interproduct_u}
\end{equation*}
and the angle between the sample eigenvector and the corresponding population subspace $\mathbb{S}_k$ in~\eqref{subspace:S} satisfies
\begin{eqnarray}
\left(\mbox{cos}\left[\mbox{angle}\left(\hat{u}_j, \mathbb{S}_k\right)\right]\right)^2=\sum_{l\in H_k} \hat{u}^2_{l,j}, \quad k=1,\cdots, r+1.
\label{subspace_u}
\end{eqnarray}

The population eigenvalues are grouped into $r+1$ tiers and
$H_k$ in \eqref{index} is the index set of the eigenvalues in the $k$th tier. Define
\begin{equation*}
\hat{U}_{k,l}=(\hat{u}_{i,j})_{i\in H_k, j\in H_l}, \quad 1\leq k, l \leq r+1.
\label{Uij}
\end{equation*}
Then, the sample eigenvector matrix $\hat{U}$
can be expressed as:
\begin{equation}
\hat{U}=[\hat{u}_1,\hat{u}_2,\cdots,\hat{u}_d]=\begin{pmatrix}
  \hat{U}_{1,1} & \hat{U}_{1,2} & \cdots & \hat{U}_{1,r+1} \\
    \hat{U}_{2,1} & \hat{U}_{2,2} & \cdots & \hat{U}_{2,r+1} \\
   \vdots    & \vdots & & \vdots  \\
  \hat{U}_{r+1,1} & \hat{U}_{r+1,2} & \cdots & \hat{U}_{r+1,r+1}
 \end{pmatrix}.
 \label{eigenvector-matrix}
\end{equation}

The proof of the asymptotic properties of the sample eigenvectors~\eqref{eigenvector:indistinct} depends on the
asymptotic properties of the sample eigenvalues, as stated in \eqref{eigenvalue:indistinct} of Theorem~\ref{Th:multiple-tiered},
which are derived in the supplementary material~\cite{shen2013online}.
The following proof considers two groups of sample eigenvectors separately.  Section~\ref{step1} obtains the
asymptotic properties for the sample eigenvectors whose index is greater than $m$.
Section~\ref{step2} derives asymptotic properties for the sample eigenvectors whose index
is less than or equal to $m$.

\subsection{Asymptotic properties of the sample eigenvectors $\hat{u}_j$ with $j>m$}\label{step1}

We derive the asymptotic properties through the following  two steps:

 \begin{itemize}
\item First, we show that as $n\rightarrow \infty$, the angle between $\hat{u}_j$ and
$u_j$ converges to 90 degrees:
\begin{equation}
\mid<\hat{u}_j,u_j>\mid^2=\hat{u}^2_{j,j}={\rm O}_{\rm a.s}\left(\frac{n}{d}\right),
\quad j=m+1,\cdots,[n\wedge d].
\label{inner:stronginconsistency}
\end{equation}

\item Then, we show that as $n\rightarrow \infty$, the angle between
$\hat{u}_j$ and the corresponding subspace $\mathbb{S}_{r+1}$ converges to 0, where $\mathbb{S}_{r+1}$ is defined as in~\eqref{subspace:S}:
\begin{equation}
\mbox{angle}<\hat{u}_j, \mathbb{S}_{r+1}>\xrightarrow{a.s} 0, \quad j=m+1,\cdots,[n\wedge d].
\label{eingvector_angle:>m}
\end{equation}
\end{itemize}

We now provide the proof for the first step.
Denote $W=\Lambda^{-\frac{1}{2}}\hat{U}\hat{\Lambda}^{\frac{1}{2}}$, where
$\hat{U}$ is the sample eigenvector matrix and $\hat{\Lambda}$ is the sample eigenvalue matrix
defined in~\eqref{sample:covariance}. It follows from~\eqref{def:z} 
 and \eqref{sample:covariance} that $WW^T=\frac{1}{n}ZZ^T$, where $Z$ is defined in~\eqref{def:z}.
Considering the $k$-th diagonal entry of the two equivalent matrices $WW^T$ and $\frac{1}{n}ZZ^T$,
and noting that
$w_{k,j}=\lambda^{-\frac{1}{2}}_k \hat{\lambda}^{\frac{1}{2}}_j \hat{u}_{k,j}$,
it follows that
\begin{equation}
\lambda^{-1}_k
\sum_{j=1}^{d} \hat{\lambda}_j\hat{u}^2_{k,j}=\sum_{j=1}^{d} w^2_{k,j} =\frac{1}{n}\sum_{i=1}^n z^2_{i,k}.
\quad k=1,\cdots,d.
\label{diagnal:kkkkk}
\end{equation}
In addition, note
 that $\frac{1}{n}\sum_{i=1}^n z^2_{i,k}\xrightarrow{\rm a.s} 1$, as $n\rightarrow\infty$,
 and
$\hat{\lambda}_j=0$ for $j>[n\wedge d]$. Combining the above with~\eqref{diagnal:kkkkk}, we obtain that
\begin{equation}
\sum_{l=1}^{r}\sum_{j\in H_l} \lambda^{-1}_k \hat{\lambda}_j\hat{u}^2_{k,j}+\sum_{j=m+1}^{[n\wedge d]} \lambda^{-1}_k\hat{\lambda}_j\hat{u}^2_{k,j} \xrightarrow{a.s} 1, \quad k=1,\cdots, d.
\label{eigenvector:value}
\end{equation}
Furthermore, it follows from~\eqref{eigenvector:value} that as $n\rightarrow \infty$,
\begin{equation}
\hat{u}^2_{j,j} \stackrel{\rm a.s}{\leq}\frac{ \lambda_j}{\hat{\lambda}_j}, \quad j=m+1,\cdots, [n\wedge d],
\label{eigenvector:jj>m}
\end{equation}
which, together with the asymptotic  properties of the sample eigenvalues~\eqref{eigenvalue:indistinct},  yields~\eqref{inner:stronginconsistency}.

We then move on to prove the second step.
According to~\eqref{subspace_u}, we need to show that
\begin{equation}
\sum_{k=m+1}^{d} \hat{u}^2_{k,j}\xrightarrow{a.s} 1, \quad j=m+1,\cdots,[n\wedge d].
\label{eigenvectors:consistency:>m}
\end{equation}
The non-zero $k$-th diagonal entry of $W^T W$ is between
its smallest and largest eigenvalues.
Since $W^T W$ shares the same non-zero eigenvalues as $\frac{1}{n}Z^TZ$, it follows that
for $j=1,\cdots, [n\wedge d]$,
\begin{equation}
\lambda_{\mbox{min}}(\frac{1}{n}Z^TZ)\leq \hat{\lambda}_j
\sum_{k=1}^{d} \lambda^{-1}_k\hat{u}^2_{k,j}=\sum_{k=1}^{d} w^2_{k,j}\leq
\lambda_{\mbox{max}}(\frac{1}{n}Z^TZ),
\label{diagkk}
\end{equation}
which yields that, for $j=m+1,\cdots,[n\wedge d]$,
\begin{equation}
\frac{\lambda_j}{\hat{\lambda}_j}\lambda_{\mbox{min}}(\frac{1}{n}Z^TZ)\leq
\sum_{k=1}^{d} \lambda_j\lambda^{-1}_k\hat{u}^2_{k,j}\leq \frac{\lambda_j}{\hat{\lambda}_j}
\lambda_{\mbox{max}}(\frac{1}{n}Z^TZ).
\label{eigenvector:>m}
\end{equation}
According to Lemma D.1 in the supplementary material~\cite{shen2013online} and
the asymptotic properties of the sample eigenvalues~\eqref{eigenvalue:indistinct},
we have that,  for $j=m+1,\cdots,[n\wedge d]$,
\begin{equation}
\frac{\lambda_j}{\hat{\lambda}_j}\lambda_{\mbox{min}}\left(\frac{1}{n}Z^TZ\right)\quad \mbox{and} \quad
\frac{\lambda_j}{\hat{\lambda}_j}
\lambda_{\mbox{max}}\left(\frac{1}{n}Z^TZ\right) \quad \xrightarrow{a.s} 1.
\label{min_max:eigenvalue:>m}
\end{equation}
In addition, it follows from Assumption $\mathcal{B}2$ that, for $j=m+1,\cdots,[n\wedge d]$,
\begin{equation}
\left\{ \begin{array}{ll}
\lambda_j\lambda^{-1}_k\rightarrow 0, \quad k=1,\cdots,m,\\
 \lambda_j\lambda^{-1}_k \rightarrow 1, \quad k=m+1,\cdots d.
\end{array} \right.
\label{eigen:>=m+1}
\end{equation}
Combining~\eqref{eigenvector:>m},~\eqref{min_max:eigenvalue:>m}, and~\eqref{eigen:>=m+1}, we have
\eqref{eigenvectors:consistency:>m}, which further leads to~\eqref{eingvector_angle:>m}.

\subsection{Asymptotic properties of the sample eigenvectors $\hat{u}_j$ with $j\in[1,m]$}\label{step2}
We need to prove that, for $j=1,\cdots, m$,
the angle between the sample eigenvector $\hat{u}_j$
 and  the corresponding population subspace
 $\mathbb{S}_l$, $j\in H_l$,
 converges to $\arccos(\frac{1}{\sqrt{1+c_l}})$,
 $l=1,\cdots, r$.
According to~\eqref{subspace_u}, we only need to show that
\begin{equation}
\sum_{k\in H_l} \hat{u}^2_{k,j}\xrightarrow{a.s} \frac{1}{1+c_l}, \quad j\in H_l, \; l=1,\cdots, r.
\label{eigenvectors_inner:<m}
\end{equation}
Below, we provide the detailed proof of~\eqref{eigenvectors_inner:<m} for $l=1$,
and briefly illustrate how  repeating the same procedure can lead to~\eqref{eigenvectors_inner:<m} for $l>2$.

In order to show~\eqref{eigenvectors_inner:<m} for $l=1$,
we need the following lemma about the
asymptotic properties of the eigenvector matrix $\hat{U}$ in~\eqref{eigenvector-matrix}:

\begin{lemma} Under Assumptions in Theorem~\ref{Th:multiple-tiered} and as $n\rightarrow\infty$,
the rows of the eigenvector matrix $\hat{U}$ satisfy
      \begin{equation}
      \sum_{l=1}^{r}(1+c_l)c_h c^{-1}_l\sum_{j\in H_l}\hat{u}^2_{k,j}\xrightarrow{a.s}1, \quad k\in H_h,\; h=1,\cdots, r,
      \label{diagk:non-noise1}
      \end{equation}
and the columns of the eigenvector matrix $\hat{U}$ satisfy
 \begin{equation}
\sum_{h=1}^r\sum_{k\in H_h}\hat{u}^2_{k,j}\xrightarrow{a.s}\frac{1}{1+c_l},
\quad j\in H_l, \; l=1,\cdots, r.
\label{eigenvector:<mm}
\end{equation}

In addition, we also have
 \begin{equation}
  \sum_{l=1}^{r}(1+c_l)\sum_{j\in H_l}\hat{u}^2_{k,j}\xrightarrow{a.s}1, \quad k\in H_1.
 \label{diagkk*}
 \end{equation}

\label{lemmaTheorem32:proof}

\end{lemma}

%
%
%
%

Lemma~\ref{lemmaTheorem32:proof} is  proven in Section D.3.3 of the supplementary material~\cite{shen2013online}. We now show how to use Lemma~\ref{lemmaTheorem32:proof} to prove
~\eqref{eigenvectors_inner:<m} for $l=1$. Let $h=1$ in~\eqref{diagk:non-noise1}, and then we have
that
\begin{equation}
\sum_{l=1}^{r}(1+c_l)c_1 c^{-1}_l\sum_{j\in H_l}\hat{u}^2_{k,j}\xrightarrow{a.s}1, \quad k\in H_1.
\label{diagk:non-noise1*}
\end{equation}
Note that $c_1 c^{-1}_l<1$ for $l>1$, and comparing~\eqref{diagkk*} with~\eqref{diagk:non-noise1*},
we get that
\begin{equation}
\sum_{l=2}^{r}\sum_{j\in H_l}\hat{u}^2_{k,j}\xrightarrow{a.s}0, \quad
 \sum_{j\in H_1}\hat{u}^2_{k,j}\xrightarrow{a.s}\frac{1}{1+c_1}, \quad k\in H_1,
\label{diag11*}
\end{equation}
which then yields that
\begin{equation}
\sum_{k\in H_1}\sum_{j\in H_1}\hat{u}^2_{k,j}\xrightarrow{a.s}\frac{q_1}{1+c_1},
\label{diag111*}
\end{equation}
where $q_1$ is the number  of eigenvalues in $H_1$~\eqref{index}.
Summing over $j\in H_1$ in~\eqref{eigenvector:<mm}, we have that
\begin{equation}
\sum_{h=1}^r\sum_{k\in H_h}\sum_{j\in H_1}\hat{u}^2_{k,j}\xrightarrow{a.s}\frac{q_1}{1+c_1}.
\label{eigenvector:<mmm*}
\end{equation}
It follows from~\eqref{diag111*} and~\eqref{eigenvector:<mmm*} that
\begin{equation}
\sum_{h=2}^r\sum_{k\in H_h}\sum_{j\in H_1}\hat{u}^2_{k,j}\xrightarrow{a.s}0,
\label{disappear*}
\end{equation}
which, together with~\eqref{eigenvector:<mm} for $l=1$, yields
\begin{equation*}
\sum_{k\in H_1} \hat{u}^2_{k,j}\xrightarrow{a.s} \frac{1}{1+c_1}, \quad j\in H_1.
\label{eigenvectors_inner:1}
\end{equation*}
which is~\eqref{eigenvectors_inner:<m} for $l=1$.

We now prove~\eqref{eigenvectors_inner:<m} for $l=2,\cdots,r$.
Note that

\begin{itemize}
\item {{it follows}} from~\eqref{disappear*} that~\eqref{diagk:non-noise1} becomes
      \begin{equation}
      \sum_{l=2}^{r}(1+c_l)c_h c^{-1}_l\sum_{j\in H_l}\hat{u}^2_{k,j}\xrightarrow{a.s}1, \quad k\in H_h,\; h=2,\cdots, r.
      \label{diagk:non-noise2}
      \end{equation}

\item {{it follows}} from~\eqref{diag11*} that~\eqref{eigenvector:<mm} becomes
 \begin{equation}
\sum_{h=2}^r\sum_{k\in H_h}\hat{u}^2_{k,j}\xrightarrow{a.s}\frac{1}{1+c_l},
\quad j\in H_l, \; l=2,\cdots, r.
\label{eigenvector:<mm2}
\end{equation}

\item similar to~\eqref{diagkk*}, we have
 \begin{equation}
  \sum_{l=2}^{r}(1+c_l)\sum_{j\in H_l}\hat{u}^2_{k,j}\xrightarrow{a.s}1, \quad k\in H_2.
 \label{diagkk2*}
 \end{equation}

\end{itemize}
Finally, combining~\eqref{diagk:non-noise2},~\eqref{eigenvector:<mm2} and~\eqref{diagkk2*}, we can prove~\eqref{eigenvectors_inner:<m} for $l=2$.
We can repeat the same procedure  for $l=3,\cdots,r$.

\begin{supplement}
\stitle{Simulations and proofs} 
\sdescription{{{The supplementary material contains additional simulation results that empirically verify the theoretical convergence of the angles between sample eigenvectors and their popularion counterparts, reported in our theorems. We also provide detailed proofs for our theorems and their extensions under both the growing sample size
 and HDLSS contexts.}}}
\slink[url]{http://www.unc.edu/$\sim$dshen/BBPCA/BBPCASupplement.pdf}
\end{supplement}

\bibliographystyle{imsart-number}
\bibliography{BBPCA}

\end{document}